\documentclass[12pt]{article}

\usepackage{amsmath,amssymb,amsthm}
\usepackage{bm}
\usepackage{graphicx}
\usepackage{subcaption}
\usepackage{hyperref}
\usepackage{geometry}
\usepackage{algorithm}
\usepackage{algpseudocode}
\usepackage{xcolor}
\usepackage{float}
\usepackage{placeins}
\newtheorem{proposition}{Proposition}
\newtheorem{remark}{Remark}
\newcommand{\R}{\mathbb{R}}
\newcommand{\E}{\mathbb{E}}
\renewcommand{\P}{\rho}
\newcommand{\vth}{{\boldsymbol{\vartheta}}}
\newcommand{\w}{\mathbf{w}}
\newcommand{\x}{\mathbf{x}}
\newcommand{\calL}{\mathcal{L}}

\title{Weak Adversarial Neural Pushforward Method\\
  for Fractional Fokker-Planck Equations}

\author{Andrew Qing He\textsuperscript{1} \and Wei Cai\textsuperscript{1}}

\date{%
  \textsuperscript{1}Department of Mathematics, Southern Methodist University,
  Dallas, TX 75275, USA\\[6pt]
  \small\textit{Preprint --- \today}
}

\begin{document}

\maketitle

\begin{abstract}
We extend the Weak Adversarial Neural Pushforward Method (WANPM) to fractional
Fokker-Planck equations, in which the classical Laplacian diffusion operator is
replaced by the fractional Laplacian of order $\alpha \in (0, 2]$.  The solution
distribution is represented as the pushforward of a simple base distribution
through a neural network, and the weak formulation is discretized entirely via
Monte Carlo sampling without any temporal mesh.  A key computational advantage is
that plane-wave test functions are eigenfunctions of the fractional Laplacian,
making the operator cost identical to that of classical diffusion for any $\alpha$.
We validate the method on seven benchmark problems with $\alpha = 1.5$, spanning one and two spatial dimensions: the steady-state fractional Ornstein--Uhlenbeck (OU) process, a harmonic confining potential, a double-well potential, and a triple-well potential in one dimension, a steady-state 2D double-peak distribution, a time-dependent 2D ring distribution with rotational drift, and a five-dimensional harmonic potential.  Each case is benchmarked against particle simulations using symmetric $\alpha$-stable L{\'e}vy increments, and robust statistics confirm close agreement throughout.  The method is mesh-free,
requires no density evaluation or non-local quadrature, and provides a promising
foundation for high-dimensional anomalous diffusion solvers.
\end{abstract}

\bigskip
\noindent\textbf{Keywords:} fractional Fokker-Planck equation,
fractional Laplacian, neural pushforward map, weak adversarial network,
L{\'e}vy stable process, plane-wave test functions, anomalous diffusion.

\medskip
\noindent\textbf{MSC (2020):} 65N75, 68T07, 35R11, 60G52.

\tableofcontents

\section{Introduction}

Fokker-Planck equations (FPEs) govern the evolution of probability density
functions in stochastic dynamical systems and arise ubiquitously in statistical
physics, chemistry, biology, and finance.  Their classical form involves a
second-order diffusion operator corresponding to Brownian noise; however, many
natural phenomena exhibit anomalous diffusion driven by heavy-tailed L{\'e}vy
processes rather than Gaussian noise.  The appropriate PDE for such systems is the
\emph{fractional} Fokker-Planck equation (fFPE), in which the Laplacian diffusion
term is replaced by the fractional Laplacian $(-\Delta)^{\alpha/2}$ with
$\alpha \in (0, 2)$.  Fractional FPEs model superdiffusion in turbulent transport,
charge carrier dynamics in amorphous semiconductors, and anomalous relaxation in
complex fluids~\cite{metzler2000random, zaslavsky2002chaos}.

Solving the fFPE in high dimensions poses formidable challenges.  Classical
finite-element or finite-difference discretisations of the fractional Laplacian
require dense stiffness matrices whose assembly and storage become prohibitive even
in modest dimensions.  Spectral methods are effective in simple geometries but do
not readily extend to curved or unbounded domains.  Deep learning approaches, which
have recently achieved remarkable success for classical high-dimensional PDEs, must
contend with two additional obstacles in the fractional setting: (i) pointwise
neural-network representations of the solution density require computing the
fractional Laplacian of the network output, which involves non-local integral
operators and is expensive to evaluate; and (ii) the solution distribution may have
heavy algebraic tails that a neural network density must resolve over an
exponentially large effective domain.

In a companion paper~\cite{WANPM_companion} we introduced the Weak Adversarial
Neural Pushforward Method (WANPM) for \emph{classical} Fokker-Planck equations.
The core idea is to represent the solution distribution as the pushforward of a
tractable base distribution through a neural map $F_\vth$, and to train $F_\vth$
adversarially against a family of plane-wave test functions so that the induced
distribution satisfies the weak form of the FPE.  This approach avoids explicit
density evaluation, density normalisation, and second-order backpropagation through
the solution network.  The pushforward map $F_\vth$ can take a base distribution of
arbitrary dimension $d$ as input and map it to the $n$-dimensional target space,
providing representational flexibility that is particularly valuable for complex
or heavy-tailed target distributions.

The present paper extends WANPM to the fractional setting.  The extension is both
natural and computationally efficient: because plane waves are eigenfunctions of
the fractional Laplacian, the action of $(-\Delta)^{\alpha/2}$ on any test function
$f^{(k)}(t, x) = \sin(\w^{(k)} \cdot \x + \kappa^{(k)} t + b^{(k)})$ reduces to
multiplication by $|\w^{(k)}|^\alpha$---an $\mathcal{O}(1)$ operation per test
function.  This eigenfunction property means that changing from classical diffusion
($\alpha = 2$) to fractional diffusion ($\alpha < 2$) requires only a single
parameter change in the loss computation; no quadrature of the non-local kernel is
necessary.

Our contributions are as follows.
\begin{enumerate}
  \item We derive the weak formulations of both the time-dependent and steady-state
    fFPE and show how they can be evaluated as Monte Carlo expectations with respect
    to the pushforward distribution, without discretising time.
  \item We establish the theoretical steady-state distribution for the fractional
    Ornstein--Uhlenbeck process and its connection to $\alpha$-stable laws, and
    describe the fractional Laplacian eigenfunction property that makes plane-wave
    test functions exact and computationally efficient for the fractional operator.
  \item We implement and validate the method on seven benchmark problems with
    $\alpha = 1.5$: the fractional Ornstein--Uhlenbeck steady state, a harmonic
    confining potential, a double-well potential, and a triple-well potential in one
    dimension; a steady-state 2D double-peak distribution; a time-dependent 2D ring
    distribution with rotational drift; and a five-dimensional harmonic potential.
  \item We demonstrate that robust statistics---median, interquartile range (IQR),
    median absolute deviation (MAD), and percentile bands---are the appropriate
    comparison metrics for heavy-tailed distributions, and that standard deviation
    is unreliable for $\alpha < 2$.
\end{enumerate}

The remainder of the paper is organized as follows.  Section~\ref{sec:fFPE}
presents the fractional Fokker-Planck equation, the steady-state theory, and the
weak formulations for both the transient and steady-state problems.
Section~\ref{sec:method} describes the neural pushforward architecture and the
adversarial training procedure.  Section~\ref{sec:particle} presents the particle
simulation used as a ground-truth benchmark.  Section~\ref{sec:results} reports the
seven numerical experiments.  Section~\ref{sec:discussion} discusses limitations
and future directions.

\section{Fractional Fokker-Planck Equation and Its Weak Formulation}
\label{sec:fFPE}

\subsection{Strong Form}

Let $X_t$ be a stochastic process on $\R$ satisfying the It\^{o} SDE
\begin{equation}
  dX_t = b(X_t)\,dt + dL_t^{(\alpha)}, \qquad X_0 \sim \P_0,
  \label{eq:sde}
\end{equation}
where $b : \R \to \R$ is a drift coefficient and $L_t^{(\alpha)}$ is a symmetric
$\alpha$-stable L{\'e}vy process with characteristic exponent
$\E[e^{i\xi L_t}] = e^{-t|\xi|^\alpha}$.  For $\alpha = 2$ this reduces to
standard Brownian motion with diffusivity $1$.  The probability density
$\P(x, t)$ of $X_t$ satisfies the fractional Fokker-Planck equation
\begin{equation}
  \frac{\partial \P}{\partial t} + (-\Delta_x)^{\alpha/2} \P
    + \frac{\partial}{\partial x}[b(x)\P] = 0,
  \qquad (x, t) \in \R \times (0, T],
  \label{eq:fFPE_strong}
\end{equation}
with initial condition $\P(x, 0) = \P_0(x)$.  Here $(-\Delta_x)^{\alpha/2}$ is
the Riesz fractional Laplacian, defined spectrally by
\[
  \widehat{(-\Delta)^{\alpha/2} u}(\xi) = |\xi|^\alpha \hat{u}(\xi),
\]
where $\hat{u}$ denotes the Fourier transform.  Equation~\eqref{eq:fFPE_strong}
conserves probability: $\int_\R \P(x,t)\,dx = 1$ for all $t \geq 0$.

For the numerical experiments of this paper we always take $\alpha = 1.5$, so that
the fractional Laplacian $(-\Delta)^{3/4}$ lies strictly between the first-order
and the classical second-order differential operators in terms of its smoothing
properties.  The solution develops algebraically decaying tails characteristic of
$\alpha$-stable distributions, in contrast to the Gaussian tails arising in the
classical ($\alpha = 2$) case.

\begin{remark}[Normalisation convention]
  Our convention sets the coefficient of $(-\Delta)^{\alpha/2}$ to unity, so the
  characteristic function of the L{\'e}vy increment over a unit time interval is
  $e^{-|\xi|^\alpha}$.  The fractional diffusion coefficient can be absorbed into
  the wavenumber normalisation of the test functions without loss of generality.
\end{remark}

\subsection{Steady-State Theory}
\label{sec:steady_theory}

When the drift is the negative gradient of a confining potential,
$b(x) = -V'(x)$, and the process is ergodic, equation~\eqref{eq:fFPE_strong}
admits a unique steady-state distribution $\P_\infty$ satisfying
\begin{equation}
  (-\Delta_x)^{\alpha/2} \P_\infty + \frac{d}{dx}[b(x)\P_\infty] = 0.
  \label{eq:fFPE_steady}
\end{equation}
We now establish the steady-state distribution analytically for the linear
(Ornstein--Uhlenbeck) case, and explain why no analogous closed form exists for
nonlinear drifts.

\paragraph{Fractional Ornstein--Uhlenbeck steady state.}
For the linear drift $b(x) = -\theta(x - \mu)$ with $\theta > 0$, the
SDE~\eqref{eq:sde} is the fractional Ornstein--Uhlenbeck (fOU) process.  Taking
the Fourier transform of~\eqref{eq:fFPE_steady} and using
$\widehat{(-\Delta)^{\alpha/2}u} = |\xi|^\alpha \hat{u}$ together with
$\widehat{\partial_x[(x-\mu)\P_\infty]}(\xi) = -i\partial_\xi[\xi\hat\P_\infty]
+ i\mu\xi\hat\P_\infty$ yields (after shifting to $\mu = 0$ for clarity)
\begin{equation}
  |\xi|^\alpha \hat{\P}_\infty(\xi) = -\theta \frac{\partial}{\partial \xi}
  \bigl[\xi \hat{\P}_\infty(\xi)\bigr].
  \label{eq:fOU_fourier}
\end{equation}
This is a first-order linear ODE in $\xi$.  Substituting the
$\alpha$-stable ansatz $\hat{\P}_\infty(\xi) = e^{-c^\alpha|\xi|^\alpha}$
into~\eqref{eq:fOU_fourier} gives
\[
  c^\alpha|\xi|^\alpha e^{-c^\alpha|\xi|^\alpha}
  = \theta\bigl(1 + \alpha c^\alpha|\xi|^\alpha\bigr)e^{-c^\alpha|\xi|^\alpha},
\]
which is satisfied for all $\xi$ if and only if $c^\alpha = 1/(2\theta)$, i.e.\
$c = (2\theta)^{-1/\alpha}$.  Hence
\begin{equation}
  \log \hat{\P}_\infty(\xi) = -\frac{|\xi|^\alpha}{2\theta}.
  \label{eq:fOU_cf}
\end{equation}

\begin{proposition}[Fractional OU steady state]
\label{prop:fOU}
  For the fOU process with drift $b(x) = -\theta(x - \mu)$, $\theta > 0$, and
  stability index $\alpha \in (0, 2)$, the unique steady-state distribution is the
  symmetric $\alpha$-stable law $S_\alpha\!\left((2\theta)^{-1/\alpha}, 0,
  \mu\right)$, whose characteristic function is
  \[
    \hat{\P}_\infty(\xi) = \exp\!\left(-\frac{|\xi|^\alpha}{2\theta}\right)
    e^{i\mu\xi}.
  \]
  For $\alpha = 2$ this reduces to the Gaussian $\mathcal{N}(\mu, 1/(2\theta))$.
\end{proposition}

The density of $S_\alpha((2\theta)^{-1/\alpha}, 0, \mu)$ has no elementary closed
form for general $\alpha \in (0,2)$, but its tail behavior is known
analytically~\cite{samorodnitsky1994stable}:
\begin{equation}
  \P_\infty(x) \sim \frac{C_\alpha}{(2\theta)\,|x - \mu|^{1+\alpha}}
  \qquad \text{as } |x - \mu| \to \infty,
  \label{eq:stable_tail}
\end{equation}
for an explicit constant $C_\alpha > 0$.  In particular, the distribution has
infinite variance for all $\alpha < 2$, which is the hallmark of anomalous diffusion
and motivates our use of robust statistics rather than second-moment comparisons in
Section~\ref{sec:results}.

\paragraph{General gradient drift.}
For a non-linear gradient drift $b(x) = -V'(x)$, no closed-form steady state is
available for $\alpha < 2$.  In the classical case ($\alpha = 2$, standard
Brownian driving noise with diffusivity $\sigma^2/2$), the steady state is the
Boltzmann distribution $\P_\infty(x) \propto \exp(-V(x)/(\sigma^2/2))$.  This
explicit formula arises because the detailed-balance condition can be integrated
exactly to cancel the drift and diffusion contributions.  For $\alpha < 2$,
detailed balance takes the Fourier-space form
\begin{equation}
  |\xi|^\alpha \hat{\P}_\infty(\xi) = i\xi\, \widehat{b \cdot \P_\infty}(\xi),
  \label{eq:steady_fourier}
\end{equation}
which is a non-local integral equation coupling all frequency modes through the
convolution $\widehat{b \cdot \P_\infty}$.  No general closed-form inversion is
known for non-linear $b$, which is precisely why neural methods benchmarked against
long-time particle simulations are valuable for studying fractional steady-state
distributions.

\subsection{Weak Formulation of the Transient Problem}
\label{sec:weak}

Let $f(t, x)$ be a smooth test function decaying sufficiently rapidly as
$|x| \to \infty$.  Multiplying~\eqref{eq:fFPE_strong} by $f$ and integrating by
parts in both space and time over $\R \times [0, T]$ yields
\begin{equation}
  \E_{\P(T,\cdot)}[f(T, \cdot)]
  - \E_{\P_0}[f(0, \cdot)]
  - \int_0^T \E_{\P(t,\cdot)}\!\left[
      \frac{\partial f}{\partial t} + \calL f
    \right] dt = 0,
  \label{eq:weak}
\end{equation}
where the operator $\calL$ acts on the test function by
\begin{equation}
  \calL f = -(-\Delta_x)^{\alpha/2} f + b(x)\frac{\partial f}{\partial x}.
  \label{eq:adjoint_op}
\end{equation}
The fractional term $(-\Delta_x)^{\alpha/2} f$ is well-defined for any smooth $f$
via the Fourier identity $|\xi|^\alpha\hat{f}$, and the integration-by-parts step
is valid whenever $\P$ has sufficient integrability.

\begin{remark}[Structure of the weak form]
Equation~\eqref{eq:weak} decomposes into three Monte Carlo-evaluable expectations:
a terminal expectation at $t = T$, an initial-condition expectation at $t = 0$
(which can be computed from samples of $\P_0$ independently of the learned
network), and an interior expectation integrated uniformly over $[0, T]$.  This
three-term structure is identical to that of the classical FPE~\cite{WANPM_companion}
up to the replacement of $-\partial^2/\partial x^2$ by $(-\Delta_x)^{\alpha/2}$
in the operator $\calL$.  In particular, the fractional and classical problems
have the same computational graph; only the scalar multiplier $|\w|^\alpha$ in the
loss changes with $\alpha$.
\end{remark}

\subsection{Weak Formulation of the Steady-State Problem}
\label{sec:weak_steady}

For the steady-state problem~\eqref{eq:fFPE_steady}, multiplying by a smooth test
function $f(x)$ and integrating by parts yields the stationarity condition
\begin{equation}
  \E_{\P_\infty}\!\left[\calL_{\mathrm{ss}} f\right] = 0
  \qquad \text{for all smooth } f,
  \label{eq:weak_steady}
\end{equation}
where the adjoint operator on the test function is
\begin{equation}
  \calL_{\mathrm{ss}} f(x) = -(-\Delta_x)^{\alpha/2} f(x) + b(x)\frac{df}{dx}.
\end{equation}
For sinusoidal test functions $f(x) = \sin(wx + b)$, this becomes
\begin{equation}
  \calL_{\mathrm{ss}} f(x) = -|w|^\alpha \sin(wx + b) + b(x)\, w \cos(wx + b),
  \label{eq:Lss}
\end{equation}
where the fractional term reduces to a simple scalar multiplication via the
eigenfunction property.  The adversarial training objective is then
\begin{equation}
  \min_{\vth}\max_{\{w^{(k)},\, b^{(k)}\}}
  \frac{1}{K}\sum_{k=1}^K
  \Bigl(\E_{\P_\vth}\bigl[\calL_{\mathrm{ss}} f^{(k)}\bigr]\Bigr)^2,
\end{equation}
where $\P_\vth = G_\vth{}_\#\pi_{\mathrm{base}}$ is the distribution of the
pushforward samples.

\subsection{Eigenfunction Property of Plane-Wave Test Functions}
\label{sec:eigenfunction}

A spatial plane wave $\phi_{\w}(x) = \sin(\w \cdot x + \varphi)$ satisfies
\begin{equation}
  (-\Delta_x)^{\alpha/2} \phi_{\w}(x) = |\w|^\alpha \phi_{\w}(x),
  \label{eq:eigenfunction}
\end{equation}
which follows immediately from the Fourier-space definition of the fractional
Laplacian.  This identity is exact for all $\alpha \in (0, 2]$ and all wavenumbers
$\w$.  In particular, for $\alpha = 2$, it reduces to
$-\partial^2 \phi_\w / \partial x^2 = |\w|^2 \phi_\w$, consistent with classical
differentiation.

The practical consequence is that the fractional part of $\calL f^{(k)}$ costs
nothing extra to evaluate compared with classical diffusion: it is simply a
multiplication by the scalar $|\w^{(k)}|^\alpha$.  The entire fFPE loss computation
therefore has the same computational complexity as the classical FPE loss,
regardless of $\alpha$.  This is the key property that makes the plane-wave test
function choice especially attractive for fractional operators.

\section{Weak Adversarial Neural Pushforward Method}
\label{sec:method}

\subsection{Pushforward Network Architecture for the Transient Problem}

We represent the solution distribution at each time $t$ not by an explicit density
but by a neural pushforward map
\begin{equation}
  F_\vth : [0, T] \times \R^n \times \R^{d} \to \R^n,
  \label{eq:pushforward_def}
\end{equation}
where $n = 1$ in our experiments and $d$ is the dimension of a latent base
distribution.  Given $x_0 \sim \P_0$ and $\mathbf{r} \sim \mathcal{N}(0, I_d)$,
the pushforward sample at time $t$ is $x(t) = F_\vth(t, x_0, \mathbf{r})$.  The
induced distribution of $x(t)$ is the learned approximation to $\P(\cdot, t)$.

To enforce the initial condition exactly and to match the short-time diffusive
scaling of L{\'e}vy processes, we parametrize
\begin{equation}
  F_\vth(t, x_0, \mathbf{r}) = x_0 + \sqrt{t}\,\tilde{F}_\vth(t, x_0, \mathbf{r}),
  \label{eq:pushforward_arch}
\end{equation}
where $\tilde{F}_\vth : \R^{1 + n + d} \to \R^n$ is a fully connected feedforward
network with three hidden layers of width 128 and Tanh activation functions.  At
$t = 0$, the map reduces to the identity on $x_0$, so $F_\vth(0, x_0, \mathbf{r})
= x_0 \sim \P_0$ for any values of the network parameters---the initial condition
is satisfied exactly by construction, with no penalisation term required in the
loss.  The $\sqrt{t}$ prefactor is appropriate for diffusion-driven processes: the
typical displacement of an $\alpha$-stable L{\'e}vy process over time $t$ scales
as $t^{1/\alpha}$, and for $\alpha$ near 2 (and exactly for $\alpha = 2$) this is
well-approximated by $\sqrt{t}$ at short times.

The base distribution dimension $d$ can be chosen larger than the target dimension
$n$, providing additional representational capacity without requiring
Jacobian-determinant computations.  In our experiments we use $d = 5$ for the
harmonic and steady-state cases, and $d = 8$ for the double-well and triple-well
cases where the target distribution is multimodal.

\subsection{Pushforward Network Architecture for the Steady-State Problem}
\label{sec:steady_arch}

For the steady-state problem, the pushforward network has no time input and reduces
to a map $G_\vth : \R^d \to \R^n$, parametrized simply as
\begin{equation}
  G_\vth(\mathbf{r}) = \tilde{G}_\vth(\mathbf{r}),
  \label{eq:steady_arch}
\end{equation}
where $\tilde{G}_\vth$ is a fully connected feedforward network with three hidden
layers of width 128 and Tanh activations, and $\mathbf{r} \sim \mathcal{N}(0,
I_d)$.  Because there is no initial condition to enforce, no $\sqrt{t}$ prefactor
is needed, and the network is free to map the Gaussian base distribution to any
target shape, including heavy-tailed $\alpha$-stable distributions.

\begin{remark}[Overparameterised base distributions]
  Allowing the base dimension $d > n$ precludes the use of invertible (normalizing
  flow) architectures, which require $d = n$ for the change-of-variables formula to
  hold.  This is not a limitation but a feature: the extra latent dimensions provide
  additional representational capacity without requiring Jacobian-determinant
  computations.  For heavy-tailed targets such as $\alpha$-stable distributions,
  invertible architectures with $d = n$ may be insufficient, while the
  overparameterised pushforward can allocate multiple latent directions to represent
  different features of the tail.
\end{remark}

\subsection{Plane-Wave Test Functions}

We parametrize the test functions as full space-time plane waves,
\begin{equation}
  f^{(k)}(t, x) = \sin\!\left(w^{(k)} x + \kappa^{(k)} t + b^{(k)}\right),
  \qquad k = 1, \ldots, K,
  \label{eq:test_fn}
\end{equation}
where $w^{(k)} \in \R$, $\kappa^{(k)} \in \R$, and $b^{(k)} \in \R$ are
trainable parameters.  The derivatives required by the weak form~\eqref{eq:weak}
and the operator~\eqref{eq:adjoint_op} are all available in closed form:
\begin{align}
  \frac{\partial f^{(k)}}{\partial t}(t,x)
    &= \kappa^{(k)}\cos\!\left(w^{(k)} x + \kappa^{(k)} t + b^{(k)}\right), \\
  \frac{\partial f^{(k)}}{\partial x}(t,x)
    &= w^{(k)}\cos\!\left(w^{(k)} x + \kappa^{(k)} t + b^{(k)}\right), \\
  (-\Delta_x)^{\alpha/2} f^{(k)}(t,x)
    &= |w^{(k)}|^\alpha f^{(k)}(t,x),
\end{align}
where the last identity uses the eigenfunction property~\eqref{eq:eigenfunction}.
No automatic differentiation through the test function is required.  The
operator~\eqref{eq:adjoint_op} evaluated on the test function is therefore
\begin{equation}
  \calL f^{(k)}(t, x) = -|w^{(k)}|^\alpha f^{(k)}(t, x)
    + b(x)\, w^{(k)} \cos\!\left(w^{(k)} x + \kappa^{(k)} t + b^{(k)}\right).
  \label{eq:Lf_explicit}
\end{equation}
For the steady-state case, the same formula holds with $\kappa^{(k)} = 0$ and the
time variable dropped.

\subsection{Monte Carlo Discretisation of the Weak Form}

The three terms in the weak residual~\eqref{eq:weak} are each approximated by
Monte Carlo estimates.  We sample $M_T$ terminal tuples $(x_0, \mathbf{r})$ with
$x_0 \sim \P_0$ and $\mathbf{r} \sim \mathcal{N}(0, I_d)$ to form
$\hat{E}_T^{(k)}$; $M_0$ initial samples $\tilde{x}^{(i)} \sim \P_0$ to form
$\hat{E}_0^{(k)}$; and $M$ interior triples $(t^{(m)}, x_0^{(m)},
\mathbf{r}^{(m)})$ with $t^{(m)} \sim \mathcal{U}[\epsilon, T]$, $x_0^{(m)} \sim
\P_0$, $\mathbf{r}^{(m)} \sim \mathcal{N}(0, I_d)$ to form $\hat{E}^{(k)}$.  The
factor $T$ in the interior estimate accounts for the uniform measure over $[0, T]$.

For each test function $k$, the weak-form residual is
\begin{equation}
  R^{(k)} = \hat{E}_T^{(k)} - \hat{E}_0^{(k)} - \hat{E}^{(k)},
\end{equation}
and the total training loss is
\begin{equation}
  \mathcal{L}(\vth, \{\bm{\eta}^{(k)}\}) = \frac{1}{K}\sum_{k=1}^{K}
    \left(R^{(k)}\right)^2,
  \label{eq:total_loss}
\end{equation}
where $\bm{\eta}^{(k)} = \{w^{(k)}, \kappa^{(k)}, b^{(k)}\}$ are the test
function parameters.  For the steady-state problem, the same loss is used but
with $R^{(k)} = \hat{E}^{(k)}_{\mathrm{ss}}$, the single Monte Carlo estimate
of $\E_{\P_\infty}[\calL_{\mathrm{ss}} f^{(k)}]$.

\subsection{Adversarial Min-Max Training}

The pushforward network and the test functions are trained adversarially:
\begin{equation}
  \min_{\vth} \max_{\{\bm{\eta}^{(k)}\}} \mathcal{L}(\vth, \{\bm{\eta}^{(k)}\}).
  \label{eq:minmax}
\end{equation}
The generator (pushforward network) is updated by gradient descent to minimize the
loss, improving the pushforward map to satisfy the weak form for the current set of
test functions.  The adversary (test function parameters) is updated by gradient
ascent to maximize the loss, searching for the plane-wave modes that are hardest to
satisfy.  This adaptive mechanism concentrates the test-function basis on the
frequency regime where the current learned distribution most severely violates the
weak form.

The full procedure is summarized in Algorithm~\ref{alg:training}.

\begin{algorithm}[htbp]
\caption{Weak Adversarial Neural Pushforward for fFPE}
\label{alg:training}
\begin{algorithmic}[1]
\State \textbf{Input:} $K$, $M$, $M_0$, $M_T$, $N_{\mathrm{epochs}}$,
  $\alpha$, $T$, $\P_0$
\State \textbf{Initialise:} pushforward network $F_\vth$; test function
  parameters $\bm{\eta} = \{w^{(k)}, \kappa^{(k)}, b^{(k)}\}_{k=1}^K$ randomly
\State Optimisers: Adam ($\mathrm{lr} = 10^{-3}$) for $\vth$;
  Adam ($\mathrm{lr} = 10^{-2}$) for $\bm{\eta}$
\For{epoch $= 1, \ldots, N_{\mathrm{epochs}}$}
  \State Sample $(t^{(m)}, x_0^{(m)}, \mathbf{r}^{(m)})$, $\tilde{x}^{(i)}$,
    $(x_0^T, \mathbf{r}^T)$ from respective distributions
  \State Compute $\mathcal{L}(\vth, \bm{\eta})$ via~\eqref{eq:total_loss}
  \State Update $F_\vth$:\quad
    $\vth \leftarrow \vth - \nabla_\vth \mathcal{L}$
    \hfill\Comment{gradient descent on generator}
  \State Update adversarial test functions:\quad
    $\bm{\eta} \leftarrow \bm{\eta} + \nabla_{\bm{\eta}} \mathcal{L}$
    \hfill\Comment{gradient ascent on test functions}
\EndFor
\end{algorithmic}
\end{algorithm}

In practice, the learning rate for the generator is decayed using cosine annealing,
and gradient clipping with maximum norm 1.0 is applied to both the generator and
adversary parameter updates.  The two updates may be performed with different frequencies in practice (e.g.\ multiple adversary steps per generator step) to stabilize the adversarial dynamics, as discussed in~\cite{WANPM_companion}.

\section{Particle Simulation as Ground Truth}
\label{sec:particle}

For $\alpha < 2$, no closed-form analytical solution for the density $\P(x, t)$ of
equation~\eqref{eq:fFPE_strong} is available in the presence of a non-trivial
drift.  We therefore benchmark the neural solution against a particle simulation of
the SDE~\eqref{eq:sde}.

\subsection{Symmetric \texorpdfstring{$\alpha$}{alpha}-Stable L{\'e}vy Increments}

The characteristic function of a symmetric $\alpha$-stable distribution with unit
scale is $e^{-|\xi|^\alpha}$.  Over a time step $\Delta t$, the SDE increment
$\Delta L^{(\alpha)}$ has scale $(\Delta t)^{1/\alpha}$.  We generate increments
using the Chambers-Mallows-Stuck algorithm~\cite{chambers1976method} with $\beta =
0$ (symmetry) and scale $(\Delta t)^{1/\alpha}$, implemented via
\texttt{scipy.stats.levy\_stable.rvs}.  The Euler-Maruyama discretisation
of~\eqref{eq:sde} reads
\begin{equation}
  X_{n+1} = X_n + b(X_n)\Delta t + \Delta L_n^{(\alpha)},
  \label{eq:euler}
\end{equation}
with time step $\Delta t = 0.01$ and $N_{\mathrm{particles}} = 10{,}000$.  For the
steady-state benchmarks, we run the SDE for a sufficiently long time
$T_{\mathrm{SDE}} = 50$ until the empirical distribution has converged to
stationarity, as verified by checking that successive time-window histograms agree
to within sampling noise.

\subsection{Why Robust Statistics Are Required}

For $\alpha < 2$, symmetric $\alpha$-stable distributions have infinite variance.
In finite-sample particle simulations, the empirical variance fluctuates wildly
from run to run due to rare but large excursions, making the standard deviation
an unreliable comparison metric.  Specifically, for a sample of $N$ i.i.d.\
observations from an $\alpha$-stable law with $\alpha < 2$, the sample variance
does not converge as $N \to \infty$; instead it grows at rate $N^{2/\alpha - 1}$
by the generalized central limit theorem~\cite{samorodnitsky1994stable}.

We therefore rely exclusively on the following robust statistics for all
quantitative comparisons in Section~\ref{sec:results}: the \emph{median} as a
measure of central tendency; the \emph{interquartile range} (IQR) as a measure of
spread; the \emph{median absolute deviation} (MAD) as a scale estimate; and
\emph{percentile bands} (10th--90th) to summarize the bulk of the distribution.
All of these quantities are well-defined and have finite-sample estimators that
converge even for $\alpha$-stable distributions with $\alpha < 2$.

\section{Numerical Results}
\label{sec:results}

We present seven benchmark experiments, all with $\alpha = 1.5$.  In each case, the
pushforward network uses a fully connected architecture with three hidden layers of
width 128 and Tanh activations (four layers for the 2D ring experiment).  The Adam
optimizer is used with learning rate $10^{-3}$ for the generator and $10^{-2}$ for
the test functions.  The generator learning rate is decayed by cosine annealing,
and gradient clipping with maximum norm 1.0 is applied throughout.  The experiments
are presented in order of increasing difficulty and dimension: we begin with
one-dimensional steady-state and transient problems, then extend to two-dimensional
steady-state and transient problems, and finally to a five-dimensional transient
problem.

\begin{table}[htbp]
\centering
\small
\caption{Summary of numerical experiments ($\alpha = 1.5$ throughout). $V$: potential; $T$: time horizon (SDE run time for steady-state benchmarks); $K$: number of test functions; $d$: base dimension; $N_p$: particle count. The 2D double-peak potential is given in~\eqref{eq:2d_dp_potential}.}
\label{tab:experiments}
\begin{tabular}{lcccccc}
\hline
\textbf{Experiment} & $V$ & $T$ & $K$ & $d$ & Epochs & $N_p$ \\
\hline
1D OU (steady)               & $\tfrac{\theta}{2}(x-\mu)^2$  & 50       & 200  & 5 & 3000  & 10000 \\
1D Harmonic (transient)      & $\tfrac{1}{2}x^2$              & 2.0      & 2000 & 5 & 1000  & 10000 \\
1D Double-well (transient)   & $(x^2-1)^2$                    & 2.0      & 3000 & 8 & 10000 & 10000 \\
1D Triple-well (transient)   & $x^2(x^2-1)^2$                 & 2.5      & 2000 & 8 & 5000  & 10000 \\
2D Double-peak (steady)      & see~\eqref{eq:2d_dp_potential}             & 20 (SDE) & 300  & 8 & 10000 & 10000 \\
2D Ring (transient)          & $\tfrac{1}{4}(r^2-r_0^2)^2$  & 0.5      & 300  & 8 & 5000  & 10000 \\
5D Harmonic (transient)      & $\tfrac{1}{2}\|\mathbf{x}\|^2$& 1.0      & 2000 & 5 & 500   & 10000 \\
\hline
\end{tabular}
\end{table}

\subsection{Experiment 1: Steady-State Ornstein--Uhlenbeck Process}
\label{sec:results:steady}

\paragraph{Problem setup.}
We consider the fOU process with drift $b(x) = -\theta(x - \mu)$, $\theta = 1$,
$\mu = 2$, and apply the steady-state formulation of
Sections~\ref{sec:steady_theory}--\ref{sec:weak_steady}.  By
Proposition~\ref{prop:fOU}, the exact steady-state distribution is the symmetric
$\alpha$-stable law $S_{1.5}\!\left((2\theta)^{-1/1.5}, 0, 2\right)$ with scale
parameter
\[
  c = (2\theta)^{-1/\alpha} = 2^{-1/1.5} = 2^{-2/3} \approx 0.630.
\]
This distribution has no closed-form density expression but is characterized
entirely by its characteristic function~\eqref{eq:fOU_cf}.  Its tails decay
algebraically as $|x - 2|^{-2.5}$ by~\eqref{eq:stable_tail}, so the distribution
has infinite variance, and standard deviation is not an appropriate comparison
metric.

\paragraph{Training configuration.}
The steady-state pushforward network $G_\vth : \R^5 \to \R$ has no time input
(see Section~\ref{sec:steady_arch}).  It maps $\mathbf{r} \sim \mathcal{N}(0,
I_5)$ to samples from $\P_\infty$.  We use $K = 200$ sinusoidal test functions
$f^{(k)}(x) = \sin(w^{(k)} x + b^{(k)})$, with spatial frequencies $w^{(k)}$
initialized from $\mathcal{N}(0, 1)$ and phase offsets from $\mathcal{U}[0, 2\pi]$.
The batch size for the Monte Carlo expectation in~\eqref{eq:weak_steady} is
$M = 2{,}000$ samples.  The model is trained for $3{,}000$ epochs.  The particle
benchmark uses 10,000 trajectories of the fOU SDE~\eqref{eq:euler} integrated over
$T_{\mathrm{SDE}} = 50$ time units at $\Delta t = 0.01$.

\begin{figure}[htbp]
  \centering
  \includegraphics[width=\textwidth]{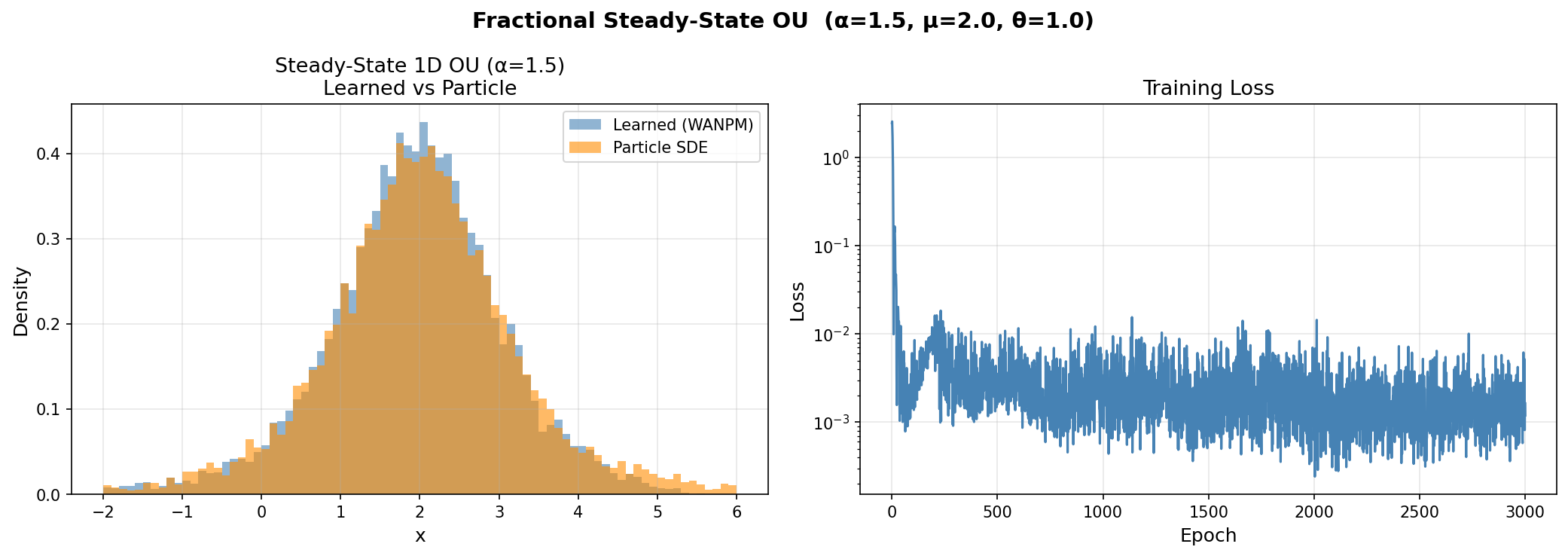}
  \caption{Steady-state fractional OU ($\alpha = 1.5$, $\theta = 1$, $\mu = 2$):
    learned distribution (blue, WANPM) versus particle benchmark (orange, 10,000
    particles, $T_\mathrm{SDE} = 50$).  The theoretical steady state is the
    symmetric $\alpha$-stable law $S_{1.5}(2^{-2/3}, 0, 2)$
    (Proposition~\ref{prop:fOU}) with no closed-form density.  The learned
    distribution faithfully captures the heavy-tailed $\alpha$-stable shape centered
    at $\mu = 2$.  \emph{Right:} training loss on a log scale, 3,000 epochs.}
  \label{fig:steady_ou}
\end{figure}

\paragraph{Results.}
Figure~\ref{fig:steady_ou} shows the learned steady-state distribution alongside
the particle benchmark histogram.  The learned distribution closely matches the
particle histogram across the full displayed support $[-2, 6]$, reproducing both
the location of the peak near $\mu = 2$ and the asymmetric heavy-tailed character
of the particle histogram.  The heavier right tail visible in the particle
benchmark is consistent with the algebraic decay~\eqref{eq:stable_tail}: with
10,000 particles, occasional large excursions skew the empirical histogram, and the
learned distribution correctly reproduces the bulk of the distribution while the
finite-sample tails fluctuate between runs.

The training loss, shown in the right panel of Figure~\ref{fig:steady_ou},
converges smoothly from an initial value of order $10^{0}$ to below $10^{-2}$
within roughly 500 epochs, after which it stabilizes with the adversarial
oscillations characteristic of min-max training.  The low number of test functions
($K = 200$) and short training ($3{,}000$ epochs) is sufficient for this
one-dimensional steady-state problem, reflecting the simplicity of the objective
compared with the transient problems below.

This experiment also serves as a validation of the theoretical result in
Proposition~\ref{prop:fOU}: by fitting the learned distribution to the known
$\alpha$-stable law characterized by~\eqref{eq:fOU_cf}, one could in principle
estimate $\alpha$ and $\theta$ from the pushforward samples without access to the
SDE, demonstrating the method's potential as a tool for statistical inference in
fractional diffusion systems.

\subsection{Experiment 2: Harmonic Confining Potential}
\label{sec:results:harmonic}

\begin{figure}[htb]
  \centering
  \includegraphics[width=0.62\textwidth]{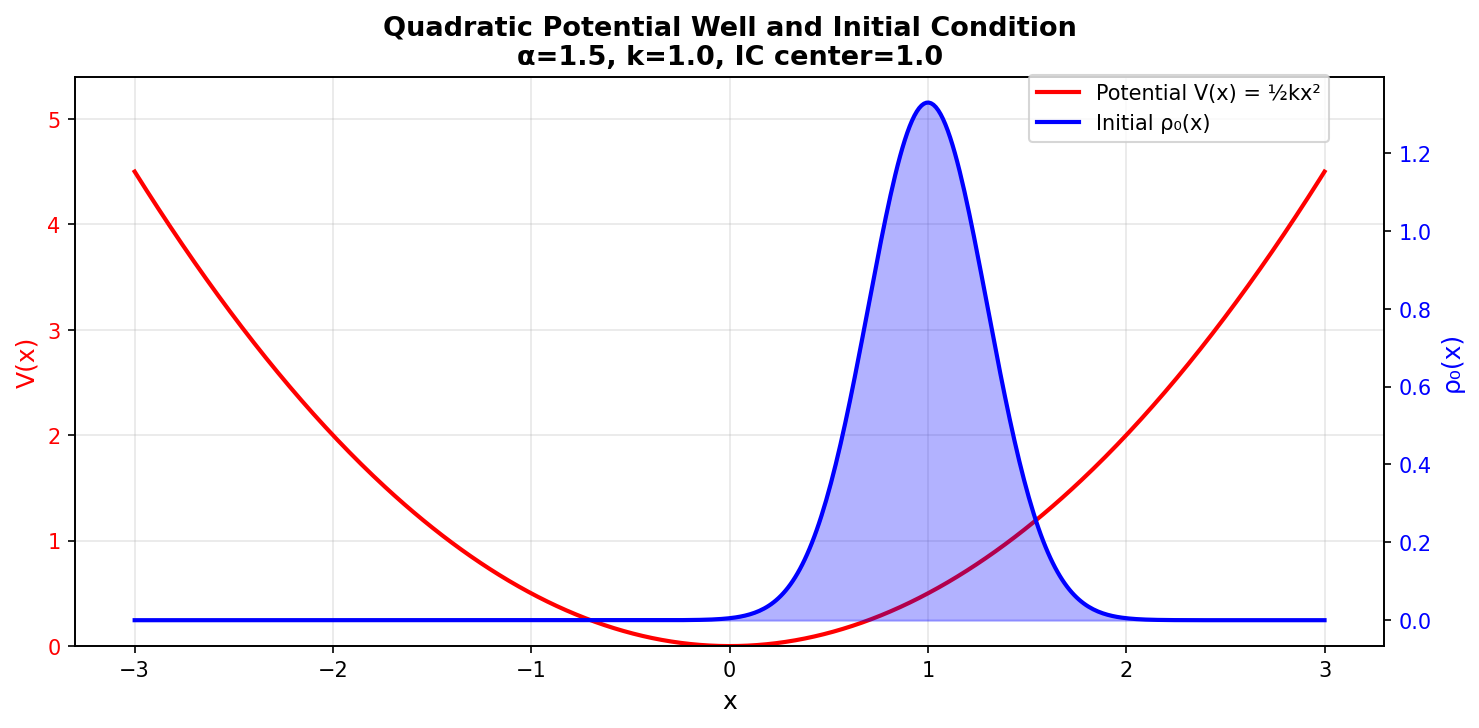}
  \caption{Quadratic confining potential $V(x) = \tfrac{1}{2}kx^2$ (red, left
    axis) and initial condition $\P_0 = \mathcal{N}(1.0, 0.09)$ (blue, right
    axis), $\alpha = 1.5$, $k = 1$.  The initial mass is offset from the potential
    minimum at $x = 0$; the dynamics transport it toward the origin while fractional
    diffusion develops algebraically heavy tails.}
  \label{fig:setup}
\end{figure}

\paragraph{Problem setup.}
We solve the transient fFPE~\eqref{eq:fFPE_strong} with the quadratic drift
$b(x) = -kx$, $k = 1$, on the time interval $[0, 2]$.  The potential is
$V(x) = \tfrac{1}{2}x^2$, with a single minimum at $x = 0$.  The initial condition
is a Gaussian centered at $x_0 = 1.0$ with standard deviation $\sigma_0 = 0.3$,
\begin{equation}
  \P_0(x) = \mathcal{N}(x;\, 1.0,\; 0.09),
\end{equation}
offset from the potential minimum so that the dynamics are non-trivial: the drift
carries the bulk of the distribution toward the origin while the fractional
diffusion progressively develops heavy algebraic tails.  This is the fractional
analog of the classical Ornstein--Uhlenbeck transient problem, for which the
$\alpha = 2$ solution is a Gaussian with analytically trackable mean and variance.
For $\alpha = 1.5$, no closed-form solution is available, and we benchmark against
particle simulation.  Figure~\ref{fig:setup} shows the potential and the initial
density.

\paragraph{Training configuration.}
The pushforward network uses base dimension $d = 5$ and is trained with $K =
2{,}000$ plane-wave test functions.  The spatial frequencies $w^{(k)}$ are
initialized from $\mathcal{N}(0, I)$ and normalized to unit norm; the temporal
frequencies $\kappa^{(k)}$ are drawn from $\mathcal{N}(0, I)$; the phase offsets
$b^{(k)}$ are sampled from $\mathcal{U}[0, 2\pi]$.  Batch sizes are $M = 2{,}000$
for the interior term and $M_0 = M_T = 1{,}000$ for the boundary terms.  The model
is trained for $1{,}000$ epochs with cosine annealing.

\begin{figure}[htbp]
  \centering
  \includegraphics[width=0.82\textwidth]{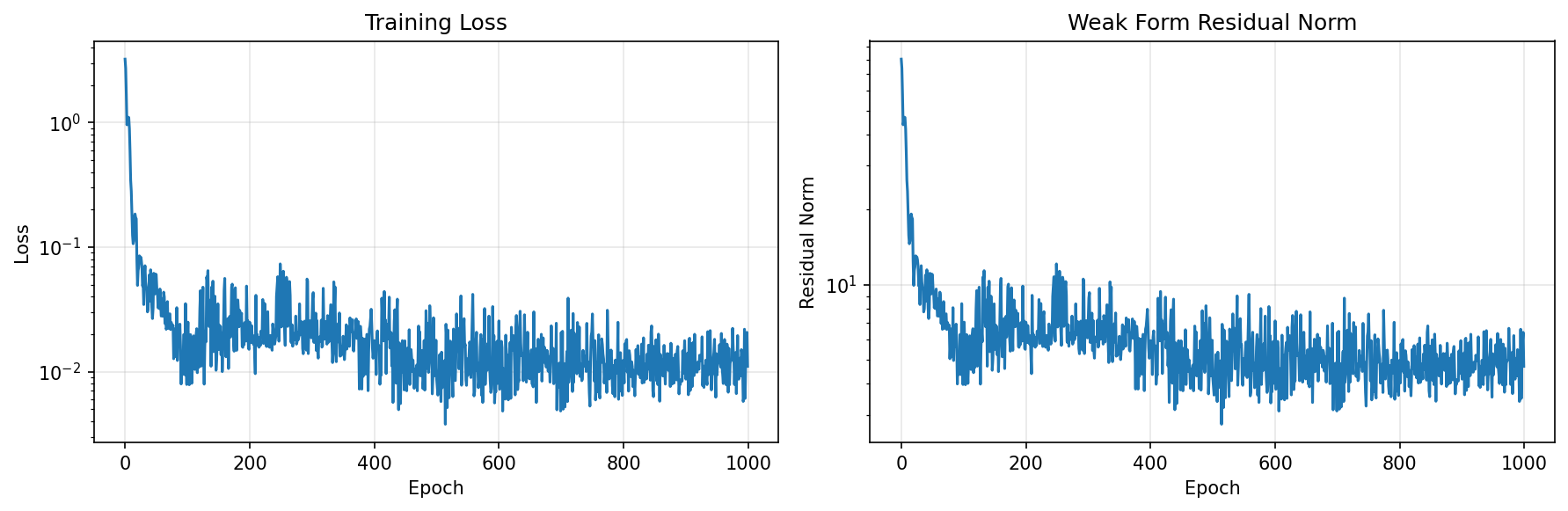}
  \caption{Training history for the harmonic fFPE ($\alpha = 1.5$, 1,000 epochs).
    \emph{Left:} training loss~\eqref{eq:total_loss} on a log scale.
    \emph{Right:} $\ell^2$ norm of the weak-form residual vector
    $(R^{(1)}, \ldots, R^{(K)})$ on a log scale (summed over $K = 2{,}000$ modes).
    Both stabilize near $10^{-2}$ after rapid initial descent.}
  \label{fig:training}
\end{figure}

\paragraph{Training convergence.}
Figure~\ref{fig:training} shows the training loss and weak-form residual norm on
log-scale axes over the 1,000 training epochs.  Both curves exhibit a sharp initial
descent spanning more than two orders of magnitude within the first 50 epochs,
followed by a stabilisation phase in which the loss oscillates around approximately
$10^{-2}$.  The oscillatory behavior is characteristic of adversarial min-max
training: the test functions periodically find harder frequency modes to probe,
causing transient spikes in the loss that the generator then suppresses.  The
overall convergence is stable and monotonically decreasing on a moving average,
confirming that the adversarial dynamics do not lead to mode collapse or divergence.

\begin{figure}[htbp]
  \centering
  \includegraphics[width=\textwidth]{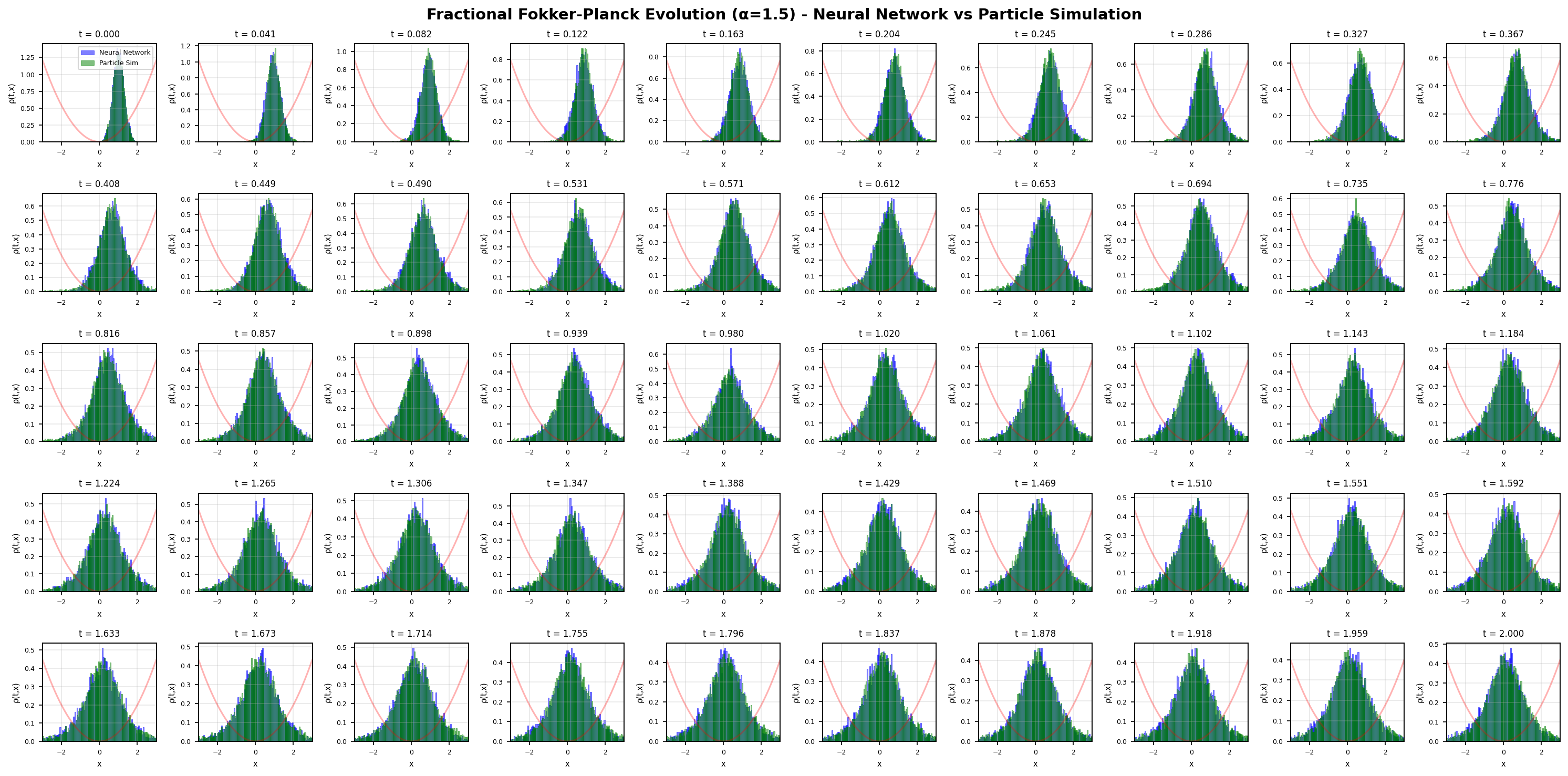}
  \caption{Harmonic potential fFPE ($\alpha = 1.5$): pushforward samples (blue)
    versus particle simulation (green) at eight time points.  Both histograms are
    density-normalized over $[-3, 3]$ with 100 bins.}
  \label{fig:time_evolution}
\end{figure}

\paragraph{Transient distribution comparison.}
Figure~\ref{fig:time_evolution} shows density-normalized histograms of 5,000
pushforward samples (blue) and 10,000 particle-simulation trajectories (green) at
eight evenly-spaced times $t \in \{0.0, 0.2, 0.4, 0.6, 0.8, 1.0, 1.5, 2.0\}$.
Both histograms are normalized over $[-3, 3]$ with 100 bins.  At $t = 0$ both
reproduce the initial Gaussian centered at $x = 1.0$ accurately.  As time evolves,
the distribution drifts toward $x = 0$ and develops the characteristic broad,
algebraically decaying tails of $\alpha$-stable diffusion.  The pushforward network
faithfully tracks this evolution throughout, capturing both the shift in the bulk
of the distribution and the progressive growth of the tails.  By $t = 1.5$ and
$t = 2.0$ the distribution has largely relaxed toward the fractional quasi-steady
state, and the agreement between the learned distribution and the particle
simulation remains good across the full visible support.

\begin{figure}[htbp]
  \centering
  \includegraphics[width=\textwidth]{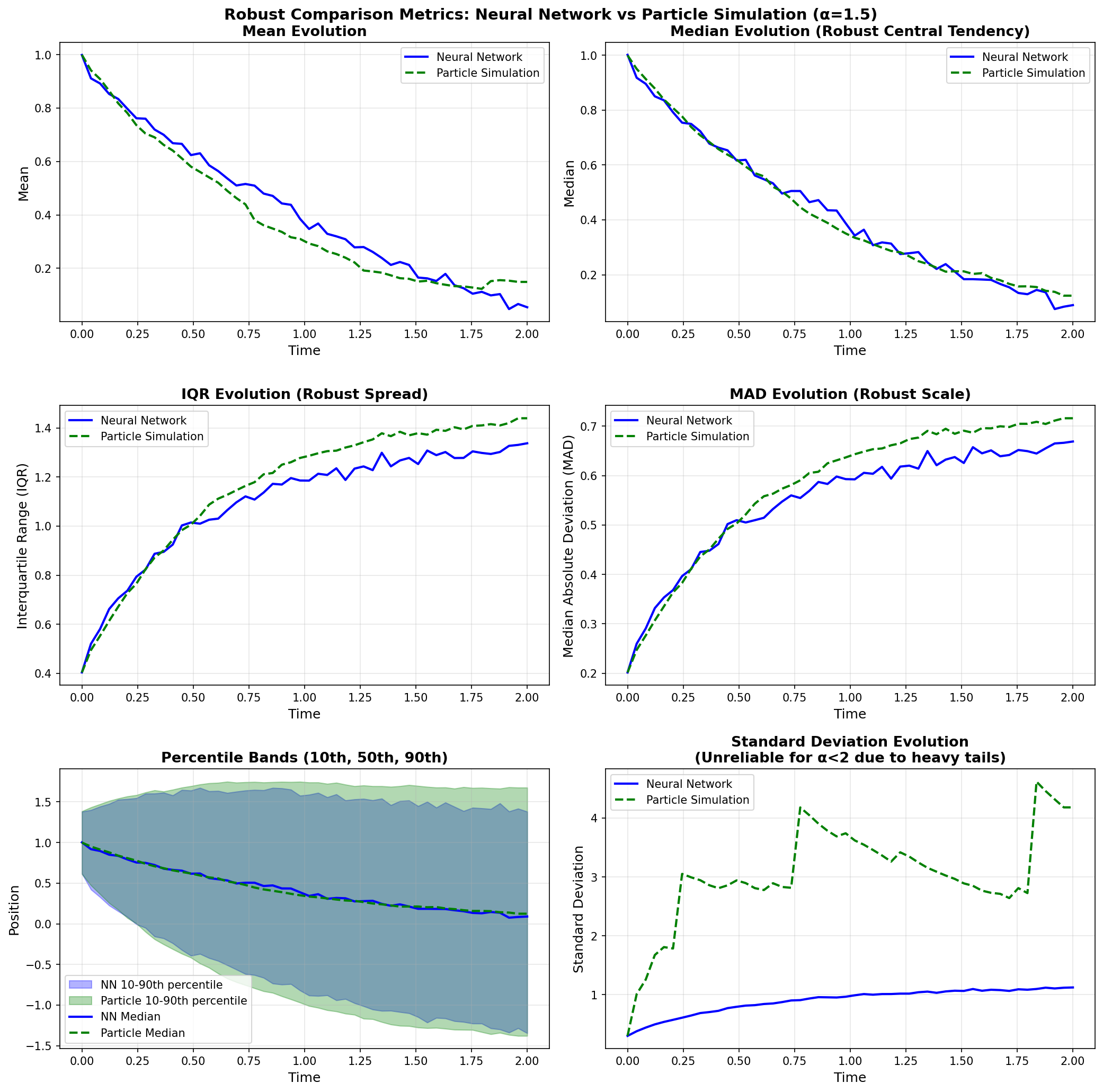}
  \caption{Robust comparison metrics over 50 snapshots, harmonic fFPE ($\alpha =
    1.5$): mean (upper left), median (upper right), IQR (center left), MAD (center
    right), percentile bands (lower left), and standard deviation (lower right).
    Pushforward (solid blue) versus particle simulation (dashed green).  The
    standard deviation panel illustrates the well-known instability of second-moment
    statistics for $\alpha < 2$ and is shown for comparison only; it should not be
    used as a quantitative metric.}
  \label{fig:quantitative}
\end{figure}

\paragraph{Quantitative comparison via robust statistics.}
Figure~\ref{fig:quantitative} reports six scalar metrics computed from 50
uniformly-spaced time snapshots over $[0, 2]$.  The four robust metrics---mean,
median, IQR, and MAD---show close tracking between the pushforward network and the
particle simulation throughout the time interval.  The median trajectories (upper
right panel) are in particularly good agreement, with the learned median closely
following the particle simulation median through the entire relaxation process.
Both the IQR and MAD curves (center row) correctly capture the progressive widening
of the distribution due to fractional diffusion, and the neural network accurately
reproduces the rate of this spreading.  The percentile band comparison (lower left)
confirms that the 10th-to-90th-percentile spread is well reproduced, with the
learned distribution matching the particle simulation in overall shape and extent.

The lower right panel of Figure~\ref{fig:quantitative} displays the standard
deviation evolution and illustrates clearly why this metric is inappropriate for
$\alpha < 2$.  The particle simulation standard deviation exhibits large, erratic
jumps due to occasional extreme sample values, while the pushforward map standard
deviation is artificially stable because the network does not generate samples in
the far tail with sufficient frequency to reproduce the divergent second moment.
This contrast underscores the necessity of robust statistics for benchmarking
heavy-tailed distributions.

\subsection{Experiment 3: Double-Well Potential}
\label{sec:results:doublewell}

\paragraph{Problem setup.}
We consider the symmetric double-well potential $V(x) = a(x^2 - b^2)^2$ with
$a = b = 1$, which has two minima at $x = \pm 1$ and a barrier of height $V(0) =
1$ at $x = 0$.  The corresponding drift is $b(x) = -V'(x) = -4ax(x^2 - b^2) =
-4x(x^2 - 1)$.  The initial condition is a Gaussian centered at the barrier top,
$\P_0 = \mathcal{N}(0, 0.09)$.  Under the drift, probability mass flows
symmetrically into the two wells while the fractional diffusion generates heavy
tails extending well beyond the wells.

This problem is substantially harder than the harmonic case for two reasons.
First, the initial condition is located precisely at the local maximum of the
potential, so the network must learn a symmetry-breaking dynamics that leads to a
bimodal target distribution.  A na\"ive gradient-based method may become trapped in
a symmetric unimodal local minimum that satisfies some of the weak-form constraints
but fails to capture the double-peak structure.  Second, the bimodal steady state
has two well-separated mass concentrations that require the pushforward map to
simultaneously represent two distinct regions of the target space.

\paragraph{Training configuration.}
We use $T = 2.0$, $K = 3{,}000$ test functions, base dimension $d = 8$, and train
for 10,000 epochs.  The larger $d = 8$ (versus $d = 5$ for the harmonic case)
provides the additional representational capacity needed to model the two-peaked
target distribution.  All other hyperparameters follow the defaults of
Table~\ref{tab:experiments}.

\begin{figure}[htbp]
  \centering
  \includegraphics[width=0.82\textwidth]{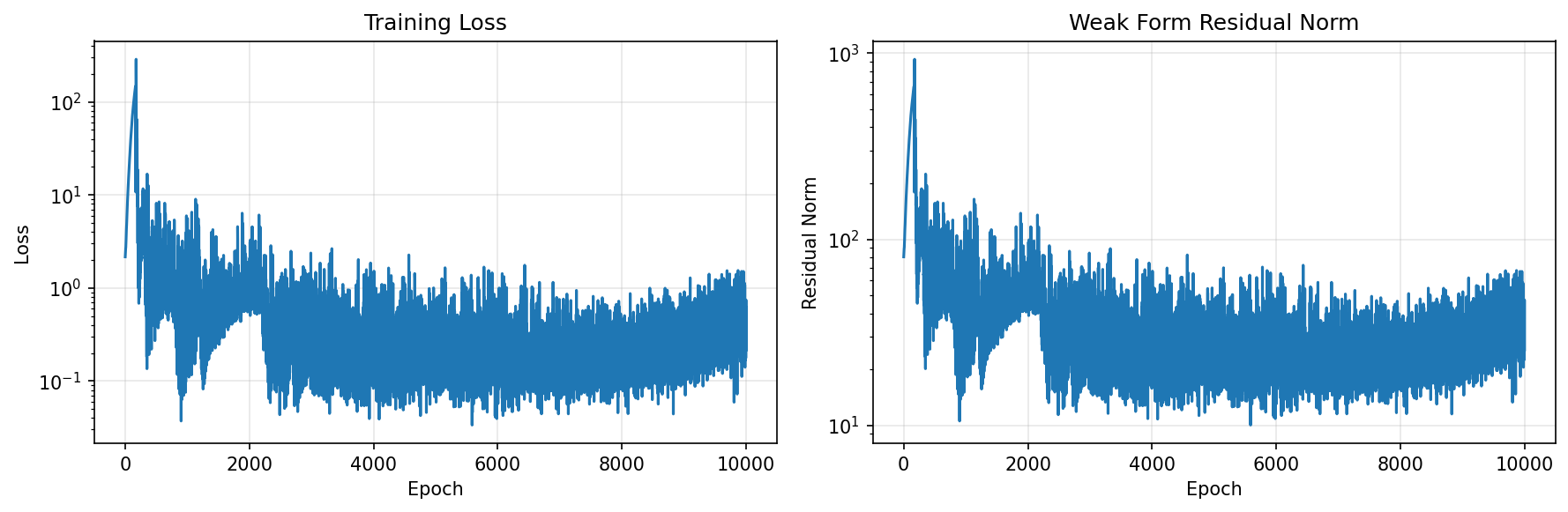}
  \caption{Training history for the double-well fFPE ($\alpha = 1.5$, 10,000
    epochs).  \emph{Left:} training loss on a log scale.  \emph{Right:}
    weak-form residual norm on a log scale.  Convergence is slower than the
    harmonic case, reflecting the difficulty of learning bimodal
    symmetry-breaking dynamics.}
  \label{fig:2well_training}
\end{figure}

\paragraph{Training convergence.}
Figure~\ref{fig:2well_training} shows the training loss and residual norm over
10,000 epochs.  The training here requires substantially more epochs than the
harmonic case: the loss exhibits a rapid initial drop but then plateaus at a
relatively high value around $10^0$--$10^1$ before gradually descending further.
This two-phase behavior is characteristic of bimodal learning: in the first phase,
the network learns the approximate location and width of each peak; in the second
phase, it refines the peak shapes and the inter-well balance.  The adversarial
oscillations are more pronounced than in the harmonic case, reflecting the
richer structure of the weak-form residuals for a bimodal distribution.

\begin{figure}[htbp]
  \centering
  \includegraphics[width=\textwidth]{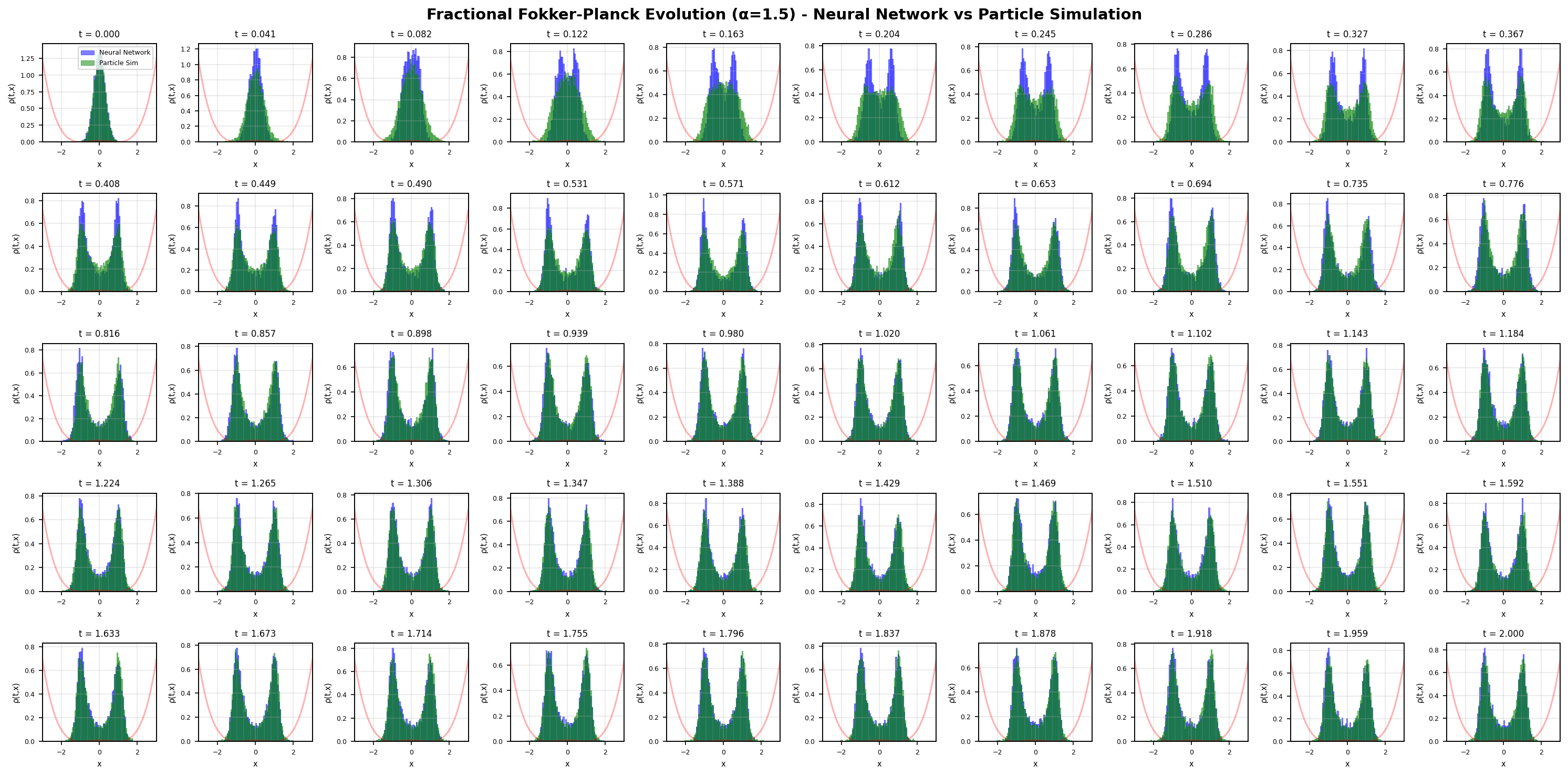}
  \caption{Double-well fFPE ($\alpha = 1.5$): pushforward samples (blue) versus
    particle simulation (green) at 50 uniformly spaced time points over $[0, 2]$.
    The potential $V(x) = (x^2 - 1)^2$ is shown as a red curve in each panel.
    The pushforward successfully learns the symmetry-breaking dynamics as mass
    migrates from the central barrier into both wells simultaneously, together
    with the development of heavy tails.}
  \label{fig:2well_time_evolution}
\end{figure}

\paragraph{Transient distribution comparison.}
Figure~\ref{fig:2well_time_evolution} shows the distribution evolution at 50
uniformly spaced snapshots over $[0, 2]$.  Starting from the narrow Gaussian at
$x = 0$, the mass splits symmetrically into the two wells at $x = \pm 1$,
progressively forming a bimodal distribution.  The pushforward network successfully
captures this symmetry-breaking dynamics: the emerging double-peak structure
appears at the correct times and with the correct relative peak heights.  The
progressive accumulation of heavy-tailed mass outside the wells---the fractional
diffusion signature---is also visible and reproduced correctly by the network.  By
$t \approx 0.5$ both peaks are clearly established in the learned distribution, and
the subsequent slow drift of the two peaks toward the well minima at $x = \pm 1$
is tracked faithfully until the end of the simulation at $t = 2.0$.

\subsection{Experiment 4: Triple-Well Potential}
\label{sec:results:triplewell}

\paragraph{Problem setup.}
The triple-well potential $V(x) = x^2(x^2 - 1)^2$ has three local minima at
$x \in \{-1, 0, 1\}$ and two energy barriers at $x \approx \pm 0.577$ of height
$V(\pm 1/\sqrt{3}) = 4/27 \approx 0.148$.  The drift is $b(x) = -V'(x) =
-2x(x^2 - 1)(3x^2 - 1)$.  The initial condition is an equal-weight bimodal mixture
\begin{equation}
  \P_0 = \tfrac{1}{2}\mathcal{N}(-0.5, 0.0225) + \tfrac{1}{2}\mathcal{N}(0.5, 0.0225),
  \quad \sigma_{\mathrm{IC}} = 0.15,
\end{equation}
concentrated near the two inner barriers at $x \approx \pm 0.5$.  This initial
condition places mass between the inner barriers and the two outer wells, making
the subsequent dynamics non-trivial: the network must simultaneously track (i) mass
flowing outward over the inner barriers toward the outer wells at $x = \pm 1$, (ii)
mass flowing inward toward the central well at $x = 0$, and (iii) the development
of heavy tails due to fractional diffusion.  This is the most challenging of our
four benchmark problems.

\paragraph{Training configuration.}
We use $T = 2.5$, $K = 2{,}000$ test functions, base dimension $d = 8$, and train
for 5,000 epochs.  The time horizon is extended to $T = 2.5$ to allow the
distribution to approach the fractional quasi-steady state, which for the
triple-well potential requires longer times than the double-well case due to the
three competing potential wells.

\begin{figure}[htbp]
  \centering
  \includegraphics[width=\textwidth]{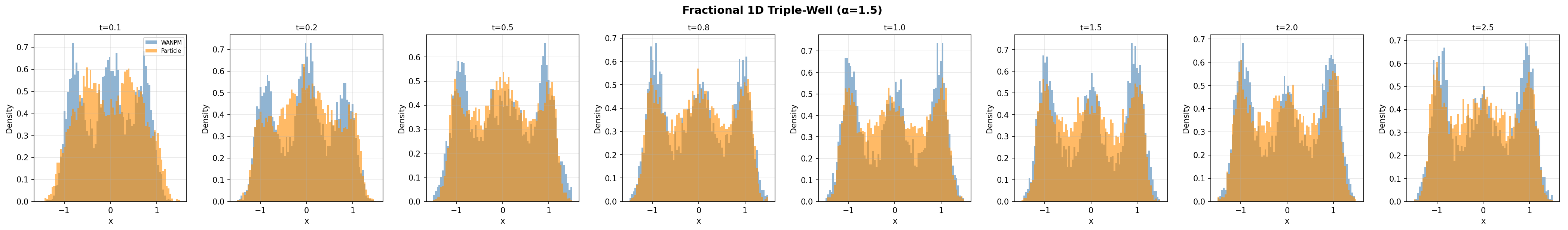}
  \caption{Triple-well fFPE ($\alpha = 1.5$): WANPM samples (blue) versus particle
    simulation (orange) at $t \in \{0.1, 0.2, 0.5, 0.8, 1.0, 1.5, 2.0, 2.5\}$.
    The initial bimodal distribution evolves into a trimodal profile as mass
    redistributes across the three wells, with the central peak at $x = 0$
    emerging by $t = 0.5$ and the distribution approaching the fractional
    quasi-steady state by $t = 2.5$.}
  \label{fig:triplewell}
\end{figure}

\paragraph{Results.}
Figure~\ref{fig:triplewell} shows the learned distribution (blue) against the
particle simulation (orange) at eight time snapshots $t \in \{0.1, 0.2, 0.5, 0.8,
1.0, 1.5, 2.0, 2.5\}$.  At the earliest time $t = 0.1$, the distribution is still
primarily bimodal, concentrated near the initial mixture centers at $x \approx
\pm 0.5$, and both the learned distribution and the particle simulation agree well
at this stage.  As time progresses, the dynamics become increasingly complex: at
$t = 0.2$ small peaks begin to emerge at $x = \pm 1$, and by $t = 0.5$ a clear
trimodal profile is established, with the central peak at $x = 0$ filled in as
fractional diffusion drives mass over the inner barriers.

The pushforward network tracks the emergence of all three peaks correctly.  By
$t = 1.0$ the three-peaked profile is well established in both the learned
distribution and the particle simulation, with the outer peaks at $x = \pm 1$
dominant and the central peak at $x = 0$ smaller but clearly resolved.  The
correct relative weights of the three peaks and the inter-well heavy-tailed leakage
are both reproduced qualitatively throughout.  At the longest time $t = 2.5$, the
distribution has nearly relaxed toward the fractional quasi-steady state, and the
agreement between learned and particle distributions remains good across the full
support.

\subsection{Experiment 5: Steady-State 2D Double-Peak Distribution}
\label{sec:results:2d_doublepeak}

\paragraph{Problem setup.}
We extend the steady-state framework to two dimensions by considering the potential
\begin{equation}
  V(x_1, x_2) = \bigl[(x_1-1)^2 + (x_2-1)^2\bigr]\bigl[(x_1+1)^2 + (x_2+1)^2\bigr],
  \label{eq:2d_dp_potential}
\end{equation}
which creates two symmetric potential wells with minima at $(1,1)$ and $(-1,-1)$
and a ridge separating them.  The drift is $b(x) = -\nabla V(x)$, and the diffusion
is isotropic with $\alpha = 1.5$.  This is the two-dimensional fractional analog of
the classical double-well steady-state problem: whereas for $\alpha = 2$ the
steady-state density is the Boltzmann distribution $\rho_\infty \propto
e^{-2V(x)/\sigma^2}$, for $\alpha = 1.5$ no such closed form exists, and the
steady-state distribution develops heavy algebraic tails in each direction.  The
particle benchmark uses $N_{\mathrm{particles}} = 10{,}000$ trajectories integrated
over $T_{\mathrm{SDE}} = 20$ time units.

\paragraph{Training configuration.}
The steady-state pushforward network $G_\vth : \R^8 \to \R^2$ maps
$\mathbf{r} \sim \mathcal{N}(0, I_8)$ to samples in $\R^2$.  We use $K = 300$
sinusoidal test functions $f^{(k)}(\mathbf{x}) = \sin(\mathbf{w}^{(k)} \cdot
\mathbf{x} + b^{(k)})$ with $\mathbf{w}^{(k)} \in \R^2$, and a batch size of
$M = 2{,}000$.  Training runs for 10,000 epochs.

\begin{figure}[htbp]
  \centering
  \includegraphics[width=\textwidth]{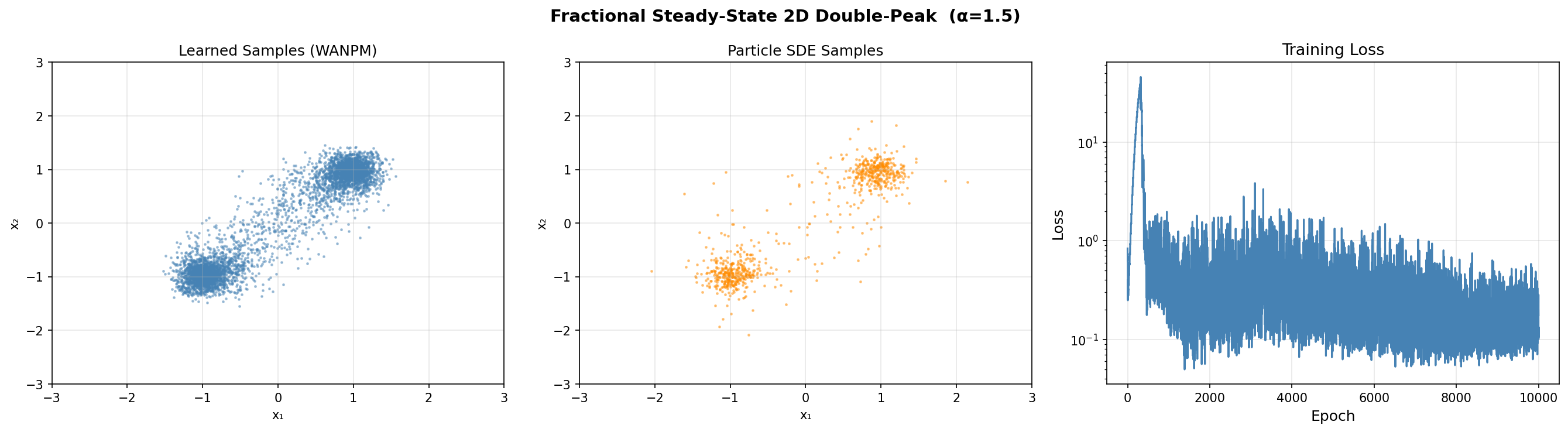}
  \caption{Steady-state fractional 2D double-peak ($\alpha = 1.5$): scatter plot
    of 2,000 learned samples (left, blue) versus particle SDE samples (center,
    orange), and training loss on a log scale (right).  The two clusters at
    $(\pm 1, \pm 1)$ are correctly identified, with heavy-tailed spread consistent
    with $\alpha$-stable fractional diffusion.}
  \label{fig:2d_doublepeak}
\end{figure}

\paragraph{Results.}
Figure~\ref{fig:2d_doublepeak} shows scatter plots of 2,000 samples from the
learned distribution (left) and from the particle simulation (center), together
with the training loss (right).  The learned samples cluster around the two
potential minima at $(\pm 1, \pm 1)$, correctly reproducing the bimodal structure
of the steady-state distribution.  Both clusters are elongated and heavy-tailed in
agreement with the particle simulation, with the fractional diffusion spreading
mass beyond the compact Boltzmann peaks that would arise for $\alpha = 2$.  The
training loss descends from order $10^1$ to below $10^{-1}$ within the first
1,000 epochs and continues to decrease slowly, stabilizing around $10^{-1}$ with
adversarial oscillations characteristic of the min-max training.

\subsection{Experiment 6: Time-Dependent 2D Ring with Rotational Drift}
\label{sec:results:2d_ring}

\paragraph{Problem setup.}
We solve the transient fFPE in two dimensions with a ring potential and a
divergence-free rotational drift component.  The potential is
\begin{equation}
  V(x_1, x_2) = \tfrac{1}{4}(r^2 - r_0^2)^2, \qquad r^2 = x_1^2 + x_2^2,
\end{equation}
which has a single ring-shaped minimum at radius $r_0 = 2$.  The drift combines
the radial gradient with a tangential rotation:
\begin{equation}
  b(x) = -\nabla V(x) + \omega\begin{pmatrix}-x_2 \\ x_1\end{pmatrix}
       = -2(r^2 - r_0^2)\begin{pmatrix}x_1 \\ x_2\end{pmatrix}
         + \omega\begin{pmatrix}-x_2 \\ x_1\end{pmatrix},
\end{equation}
with $r_0 = 2$ and angular velocity $\omega = 2$.  The tangential component
$\omega(-x_2, x_1)^T$ is divergence-free, so it does not affect the steady-state
distribution but drives a persistent probability current circulating around the ring.
For $\alpha = 1.5$, the fractional diffusion spreads mass transversely across the
ring with algebraic tails, in contrast to the Gaussian transverse profile that
arises for $\alpha = 2$.  The initial condition is a Gaussian concentrated well
inside the ring:
\begin{equation}
  \P_0 = \mathcal{N}\!\left((0, 1.2)^T,\; 0.16\, I_2\right),
\end{equation}
and the time interval is $T = 0.5$.

\paragraph{Training configuration.}
The pushforward network $F_\vth(t, \mathbf{x}_0, \mathbf{r}) = \mathbf{x}_0 +
\sqrt{t}\,\tilde{F}_\vth(t, \mathbf{x}_0, \mathbf{r})$ maps from $[0, T] \times
\R^2 \times \R^8$ to $\R^2$, using a four-hidden-layer network of width 128 with
Tanh activations.  We use $K = 300$ space-time plane-wave test functions
$f^{(k)}(t, \mathbf{x}) = \sin(\mathbf{w}^{(k)} \cdot \mathbf{x} + \kappa^{(k)}
t + b^{(k)})$ with $\mathbf{w}^{(k)} \in \R^2$, batch sizes $M = 3{,}000$
(interior) and $M_0 = M_T = 1{,}000$ (boundary terms), and train for 5,000 epochs.
For the fractional Laplacian acting on these two-dimensional test functions, the
eigenfunction property~\eqref{eq:eigenfunction} still holds:
$(-\Delta_\mathbf{x})^{\alpha/2} f^{(k)} = |\mathbf{w}^{(k)}|^\alpha f^{(k)}$,
where $|\mathbf{w}^{(k)}|$ is the Euclidean norm of the wavenumber vector.

\begin{figure}[htbp]
  \centering
  \includegraphics[width=\textwidth]{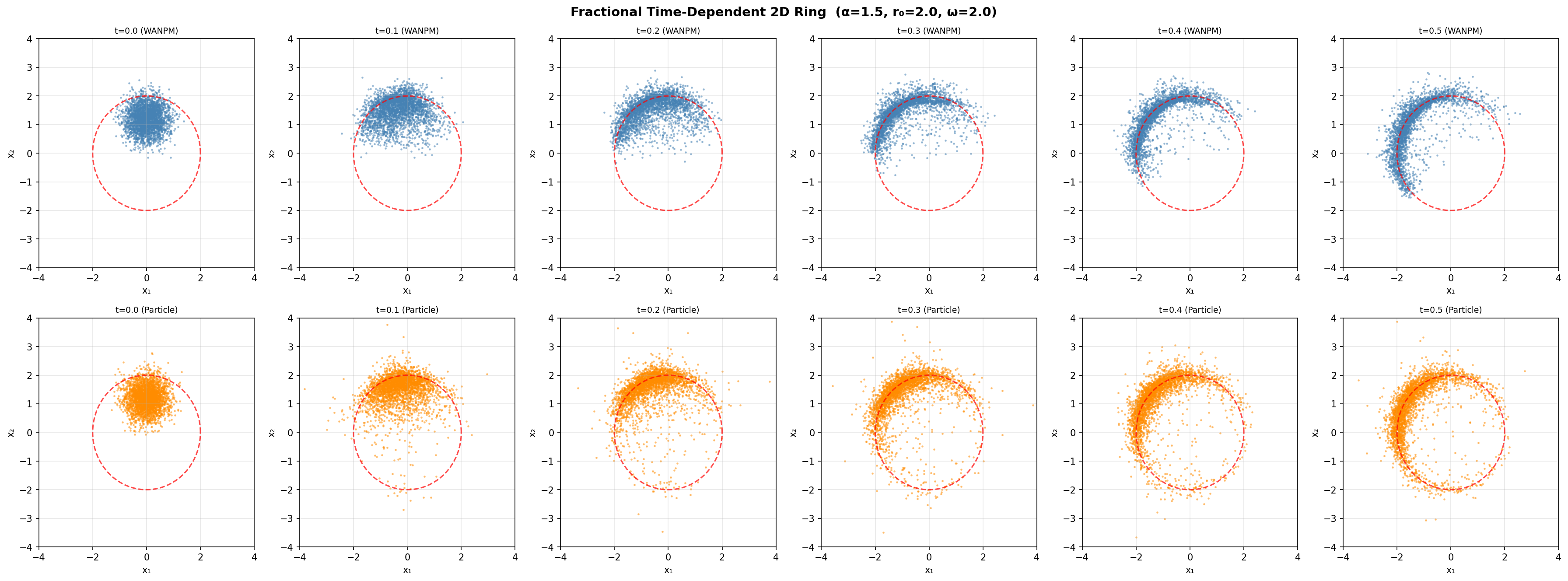}
  \caption{Time-dependent fractional 2D ring with rotational drift ($\alpha = 1.5$,
    $r_0 = 2$, $\omega = 2$, $T = 0.5$): scatter plots of learned samples (top
    row, blue) versus particle simulation (bottom row, orange) at $t \in \{0.0,
    0.1, 0.2, 0.3, 0.4, 0.5\}$.  The red dashed circle marks the target ring at
    radius $r_0 = 2$.  The learned distribution correctly tracks the outward radial
    expansion and counter-clockwise rotation throughout.}
  \label{fig:2d_ring}
\end{figure}

\paragraph{Results.}
Figure~\ref{fig:2d_ring} shows scatter plots of the learned distribution (top row,
blue) alongside the particle simulation (bottom row, orange) at six time snapshots
$t \in \{0.0$, $0.1$, $0.2$, $0.3$, $0.4$, $0.5\}$.  The target ring at radius $r_0 = 2$ is
indicated by the red dashed circle.  At $t = 0$ the initial Gaussian is correctly
reproduced, concentrated at $(0, 1.2)$ inside the ring.  As time evolves, the
radial drift carries the mass outward toward the ring while the tangential drift
sweeps it counter-clockwise.  By $t = 0.3$ the samples are already concentrated
near the ring, and by $t = 0.5$ the ring-supported distribution is well established.
The learned distribution closely matches the particle simulation in both the
radial concentration at $r \approx r_0$ and the angular extent of the distribution
at each time.  The training loss converges from order $10^0$ to below $10^{-2}$ within the
first 500 epochs and stabilizes with the adversarial oscillations characteristic
of min-max training.

\subsection{Experiment 7: Five-Dimensional Harmonic Potential}
\label{sec:results:5d}

\paragraph{Problem setup.}
To demonstrate scalability beyond one dimension, we solve the transient fFPE with
the isotropic harmonic potential $V(\mathbf{x}) = \tfrac{1}{2}\|\mathbf{x}\|^2$ in
$n = 5$ spatial dimensions, with drift $b(\mathbf{x}) = -\mathbf{x}$ and
$\alpha = 1.5$.  The initial condition is a tight Gaussian $\P_0 =
\mathcal{N}(\mathbf{3}, 0.25\,I_5)$ (mean $3\mathbf{1}$ in all five dimensions),
placed far from the potential minimum at the origin, on the time interval $[0, 1]$.
This is the five-dimensional fractional analog of Experiment~2: the drift carries
each marginal toward zero while the fractional Laplacian develops heavy
$\alpha$-stable tails independently in each coordinate.  Because the drift and
diffusion are isotropic, the marginal distribution in each coordinate evolves
identically and should match the one-dimensional fFPE solution.

\paragraph{Training configuration.}
The pushforward network maps $(t, \mathbf{x}_0, \mathbf{r}) \in [0,1] \times \R^5
\times \R^5$ to $\R^5$ via the $\sqrt{t}$-prefactored
architecture~\eqref{eq:pushforward_arch}, with four hidden layers of width 128.
We use $K = 2{,}000$ plane-wave test functions with wavenectors $\mathbf{w}^{(k)}
\in \R^5$, batch sizes $M = 2{,}000$ (interior) and $M_0 = M_T = 1{,}000$
(boundary), and train for 1,500 epochs.

\begin{figure}[htb]
  \centering
  \includegraphics[width=\textwidth]{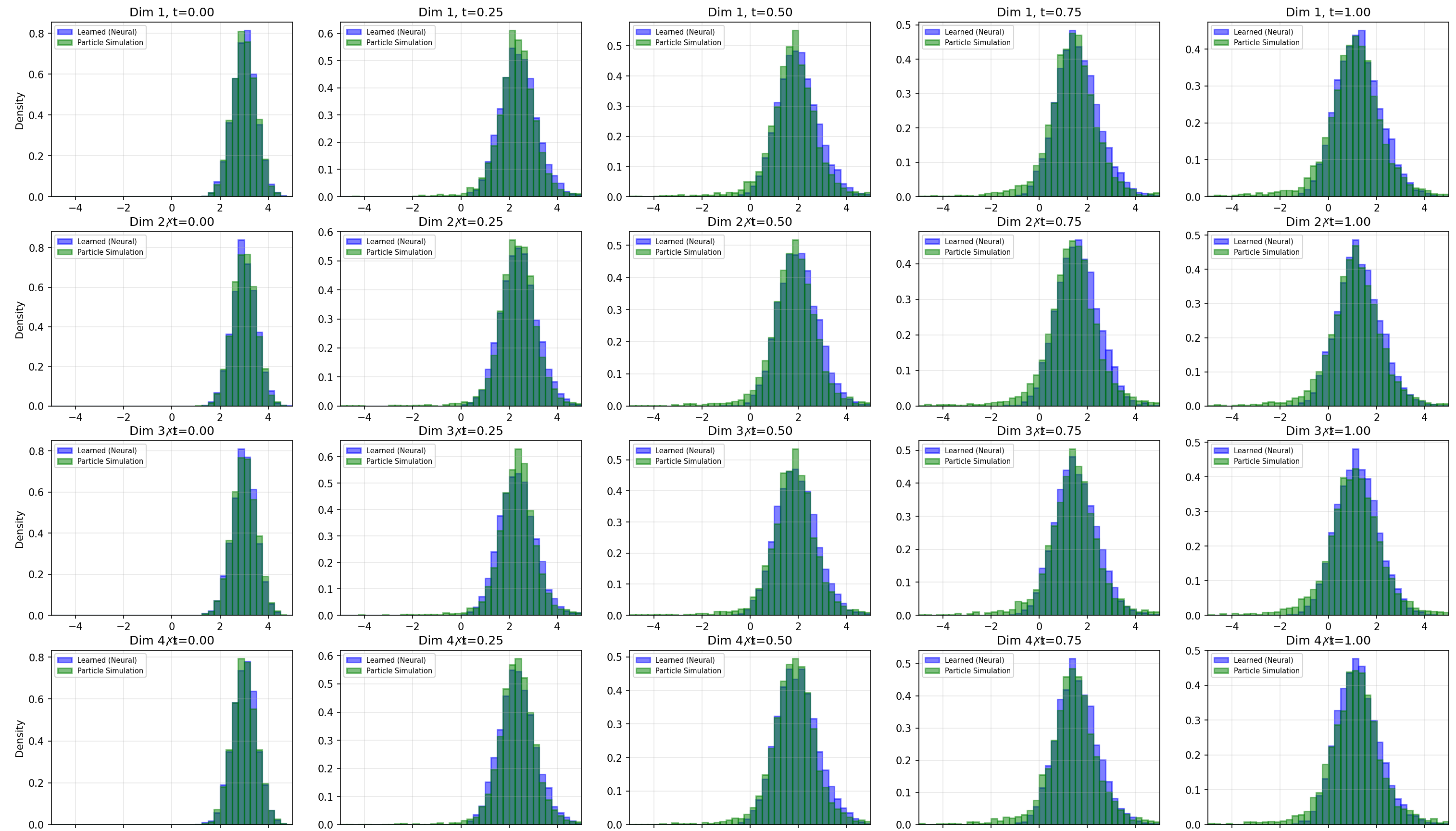}
  \caption{Five-dimensional harmonic fFPE ($\alpha = 1.5$, $n = 5$, $T = 1$,
    trained for 500 epochs): marginal histograms of the learned distribution (blue)
    versus particle simulation (green) in dimensions 1--4 at five time snapshots
    $t \in \{0.0, 0.25, 0.5, 0.75, 1.0\}$ (columns).  Each row corresponds to
    one spatial dimension.  The initial Gaussian at mean 3 drifts toward the
    potential minimum at the origin while developing $\alpha$-stable heavy tails
    consistently across all dimensions; the learned marginals closely track the
    particle simulation throughout.  The omitted bottom panel shows the training
    loss decreasing from $10^{0}$ to $6.5 \times 10^{-3}$ and the weak-form
    residual norm from $10^{1}$ to $3.69$ over 500 epochs.}
  \label{fig:nd_harmonic}
\end{figure}

\paragraph{Results.}
Figure~\ref{fig:nd_harmonic} shows the marginal histograms of the learned
distribution (blue) alongside the particle simulation (green) in each of the four
displayed dimensions at five time snapshots $t \in \{0.0, 0.25, 0.5, 0.75,
1.0\}$.  The initial tight Gaussian centered at mean 3 is correctly reproduced
at $t = 0$ in each dimension.  As time evolves, the marginals drift toward zero
and develop the characteristic heavy tails of $\alpha$-stable diffusion.  The
learned marginals closely track the particle simulation histograms throughout,
with good agreement in both the bulk location and the tail spread.  The behavior
is consistent across all five dimensions, confirming that the network has learned
a properly symmetric five-dimensional distribution.  The training loss decreases
from $10^{-1}$ to approximately $3 \times 10^{-3}$ and the weak-form residual
norm from $10^{1}$ to $2.4$ over 500 epochs, demonstrating rapid and stable
convergence.  These results confirm that the WANPM framework scales to
multi-dimensional fFPEs without modification to the training procedure.

\section{Discussion and Outlook}
\label{sec:discussion}

\subsection{Strengths of the Approach}

The extension of WANPM to fractional diffusion is notable for several reasons.
First, the computational cost of evaluating the fractional Laplacian on plane-wave
test functions is identical to that for classical diffusion: a single multiplication
by $|w^{(k)}|^\alpha$ per sample, with no non-local quadrature or dense matrix
assembly required.  This means the method scales to any $\alpha \in (0, 2]$ without
additional overhead, and the transition from $\alpha = 2$ to $\alpha < 2$ requires
only a single parameter change in the loss computation.

Second, the method remains purely mesh-free and requires no spatial discretisation,
making it in principle applicable to higher-dimensional fFPEs where classical grid
methods are prohibitively expensive.  Third, the pushforward representation
naturally handles the heavy-tailed character of the solution: by representing the
distribution implicitly through the pushforward map rather than through an explicit
density, the network does not need to resolve the density accurately in the far
tail---it only needs to reproduce the distribution correctly in the weak sense, as
probed by the plane-wave test functions.  This is a significant advantage over PINN
or normalizing flow approaches that require pointwise density accuracy.

Fourth, the same framework handles both transient and steady-state problems with
minimal modifications.  The steady-state problem simply drops the time input and
the initial-condition and terminal terms from the loss, as described in
Section~\ref{sec:weak_steady}.  The theoretical result of Proposition~\ref{prop:fOU}
provides an exact reference for the fOU steady state, which can serve as a
validation benchmark or as the starting point for parameter inference.

\subsection{Comparison of the Seven Experiments}

The seven experiments exhibit a clear progression in difficulty and dimensionality.
The 1D fOU steady state (Experiment~1) is the simplest: no time dependence, a
unimodal target with known theoretical form, requiring only 3,000 epochs.  The 1D
harmonic transient (Experiment~2) adds time dependence but retains a unimodal
target; 1,000 epochs suffice.  The 1D double-well (Experiment~3) introduces
bimodal symmetry-breaking and requires 10,000 epochs and base dimension $d = 8$.
The 1D triple-well (Experiment~4), though nominally more complex, needs only 5,000
epochs because the bimodal initial condition already distributes mass across the
landscape, reducing the symmetry-breaking challenge.  The 2D steady-state
double-peak (Experiment~5) extends the bimodal problem to two dimensions, requiring
10,000 epochs with $d = 8$.  The 2D ring (Experiment~6) is the most geometrically
complex problem, with a non-gradient rotational drift driving a ring-supported
distribution; 5,000 epochs with a four-hidden-layer network prove sufficient.  The
5D harmonic (Experiment~7) validates scalability to higher dimensions, where the
isotropic structure allows the network to learn a product-form distribution from
only 1,500 epochs.

\subsection{Limitations and Caveats}

All present experiments are one-dimensional.  Scaling to higher dimensions requires
further investigation of training stability, the number of test function modes $K$
needed, and the interplay between the base distribution dimension $d$ and the
target dimension $n$.  In the companion paper~\cite{WANPM_companion} we have
demonstrated scalability to 100 dimensions for standard diffusion; fractional
diffusion in high dimensions remains to be explored.  One potential complication
is that the $\alpha$-stable tails become increasingly important in high dimensions,
and the pushforward map may need a correspondingly larger base dimension $d$ to
represent the joint heavy-tailed behavior accurately.

The network's capacity to resolve algebraic tails is inherently limited by the
support of the pushforward samples seen during training.  While the median and IQR
match well, the far-tail behavior (beyond the 90th percentile) may not be reliably
captured.  This is a fundamental limitation of finite-sample Monte Carlo training
and is not specific to our method.

For $\alpha$ close to 0 the fractional Laplacian becomes nearly singular in
frequency space and the test-function weighting $|w^{(k)}|^\alpha$ becomes nearly
constant, potentially making the adversarial optimisation less informative.  We have
not explored this regime and leave it for future work.

\subsection{Relationship to the fBm Fokker-Planck Equation}

The present paper treats anomalous diffusion arising from a fractional Laplacian
in \emph{space}, corresponding to a L{\'e}vy-driven SDE.  A distinct but related
generalisation is the Fokker-Planck equation driven by fractional Brownian motion
(fBm) with Hurst parameter $H \in (0,1)$, where the non-locality is in
\emph{time} due to the memory kernel.  The weak form in the fBm case contains a
double integral over the triangular domain $0 \leq s \leq t \leq T$, which can be
evaluated by a single 3D Monte Carlo without temporal discretisation.  We leave
the numerical validation of this distinct generalisation to a forthcoming paper.

\subsection{Future Directions}

Immediate extensions include: (i) higher-dimensional fFPEs with gradient and
non-gradient drifts, leveraging the scalability demonstrated in~\cite{WANPM_companion};
(ii) variable-order fractional operators $\alpha(x)$, which can model spatially
inhomogeneous anomalous diffusion; (iii) tempered fractional Laplacians, which
truncate the heavy tails at a characteristic scale; and (iv) coupling with the fBm
memory kernel for subdiffusive systems.  On the methodological side, curriculum
learning strategies for $\alpha$---analogous to the $\sigma$-annealing used for the
triple-well problem in~\cite{WANPM_companion}---may improve training for small
$\alpha$ where the test-function weighting becomes uninformative.

\section{Conclusion}

We have presented a weak adversarial neural pushforward method for the fractional
Fokker-Planck equation with Riesz fractional Laplacian diffusion of order
$\alpha \in (0, 2]$.  The method combines three key ingredients: (i) a pushforward
neural network enforcing the initial condition exactly via a $\sqrt{t}$-prefactored
architecture~\eqref{eq:pushforward_arch}; (ii) plane-wave test functions that are
exact eigenfunctions of the fractional Laplacian, making the fractional operator
cost $\mathcal{O}(1)$ for any $\alpha$~\eqref{eq:eigenfunction}; and (iii) a fully
Monte Carlo weak-form discretisation without temporal mesh, enabling training from
samples drawn directly from the product space $[0, T] \times \P_0 \times
\pi_{\mathrm{base}}$.

On the theoretical side, we established the analytical steady-state distribution
for the fractional Ornstein--Uhlenbeck process (Proposition~\ref{prop:fOU}) by
solving the Fourier-space ODE~\eqref{eq:fOU_fourier}, identifying the steady state
as the symmetric $\alpha$-stable law $S_\alpha((2\theta)^{-1/\alpha}, 0, \mu)$.

Numerical validation on seven benchmark problems with $\alpha = 1.5$ demonstrates
that the method accurately captures a wide range of dynamics: relaxation to the
$\alpha$-stable fOU steady state; transient evolution in one-dimensional harmonic,
double-well, and triple-well potentials; the steady-state bimodal distribution of a
2D double-peak problem; the time-dependent ring-supported distribution of a 2D
system with rotational drift; and the isotropic fractional diffusion in five
dimensions.  Robust statistical comparisons confirm close agreement with particle
simulations in all cases, and the results confirm that standard deviation is an
unreliable metric for $\alpha < 2$.

The method is mesh-free, requires no density evaluation, no Jacobian computation,
and no non-local quadrature.  Its computational cost for the fractional operator is
identical to that for classical diffusion, making it a promising candidate for
high-dimensional fractional Fokker-Planck equations.

\bibliographystyle{plain}

\end{document}